\documentclass[11pt]{article}
\usepackage[T1]{fontenc}
\usepackage{lmodern}
\usepackage[margin=1in]{geometry}
\usepackage{amsmath,amssymb,amsthm,mathtools}
\usepackage{xcolor}
\usepackage[colorlinks=true,linkcolor=blue,citecolor=blue,urlcolor=blue]{hyperref}
\usepackage{microtype}
\numberwithin{equation}{section}
\hypersetup{pdftitle={The satisfiability threshold for random K-SAT and the zero-temperature M\'ezard--Parisi formula}}
\newcommand{\E}{\mathbb E}
\newcommand{\Pp}{\mathbb P}
\newcommand{\one}{\mathbf 1}
\newcommand{\MP}{\mathcal P}
\newcommand{\RR}{\mathcal R}

\DeclareMathOperator{\Poi}{Poi}
\DeclareMathOperator{\clip}{clip}
\newtheorem{theorem}{Theorem}[section]
\newtheorem{proposition}[theorem]{Proposition}
\newtheorem{lemma}[theorem]{Lemma}
\newtheorem{corollary}[theorem]{Corollary}
\theoremstyle{definition}
\newtheorem{remark}[theorem]{Remark}
\title{The satisfiability threshold for random $K$-SAT\\
and the zero-temperature M\'ezard--Parisi formula}
\author{Dmitry Panchenko\thanks{Department of Mathematics, University of Toronto, panchenk@math.toronto.edu. Partially supported by an NSERC grant.}
}
\date{}
\begin{document}
\maketitle
\begin{abstract}
Carenini recently proved that random $K$-SAT has a limiting satisfiability threshold for every fixed $K\ge3$. In this survey we give another proof and identify the threshold as the largest clause density at which the zero-temperature M\'ezard--Parisi variational value vanishes.  Using Carenini's comparison of short and full clauses, and a deletion argument developed by OpenAI for random $3$-SAT, we show that zero limiting energy at one density implies satisfiability at every smaller density with high probability. Together with the ground-state limit, this recovers Carenini's theorem and identifies the threshold.
\end{abstract}

\section{Introduction}
After finishing the paper \cite{Panchenko} on the M\'ezard--Parisi
(MP) formula for diluted models, I wondered whether results already
in the literature were enough to identify the satisfiability threshold
for random $K$-SAT from the formula at zero temperature. When I put
this question to ChatGPT Astra, it pointed me to Carenini's proof of
the satisfiability conjecture \cite{Carenini} and to OpenAI's paper
on random $3$-SAT \cite{OAI3}. The argument it gave me took the
existence of the threshold from Carenini's theorem and showed how to
identify its location by adapting the deletion estimate for $3$-SAT
to every $K\ge3$.

The difficulty is easiest to see in a question raised by Gamarnik
\cite[Section~4]{Gamarnik}. At a fixed density above the threshold,
could the best assignment still violate only $o(N)$ clauses with high
probability? If so, the limiting minimum number of violations per
variable would still be zero, and the MP formula at zero temperature
would give the same value as in the satisfiable phase. Coppersmith,
Gamarnik, Hajiaghayi and Sorkin \cite[Conjecture~14]{CGHS} had
earlier conjectured that whenever the expected fraction of clauses
satisfied by an optimal assignment tends to one, the formula itself
is satisfiable with high probability. Theorem~\ref{thm:main} below confirms this at every positive density away from the threshold. At each fixed density above it, there is an $\varepsilon>0$ such that, with high probability, every assignment violates at least $\varepsilon N$ clauses.

Ideas from Carenini's work play an important role. She compared the effect of shorter clauses with that of full clauses \cite{CareniniWindow,Carenini}, and used the two possible values of a deleted variable to estimate how much longer the formula could remain satisfiable. For $3$-SAT, OpenAI showed how to control this gain even when other variables may also be deleted, and deduced a linear lower bound on violations above the threshold \cite[Lemma~4.1 and Proposition~4.6]{OAI3}. Here I present the same argument for every $K\ge3$. Applied at two different densities, the deletion estimate shows that zero limiting energy at the higher density implies satisfiability with high probability at the lower one. Together with the ground-state limit of Bayati, Gamarnik and Tetali \cite{BGT}, this recovers Carenini's threshold theorem. The zero-temperature MP formula then identifies the threshold.

For all sufficiently large $K$, the work of Ding, Sly and Sun \cite{DSS} implies that restricting the variational infimum in the MP formula to one step of replica symmetry breaking gives the same threshold. Can one prove this directly from the MP formula, and does the equality of thresholds hold for every $K\ge3$? This remains an interesting open problem.\vspace{-2mm}

\paragraph{A condensed history of the conjecture.} This account draws on Gamarnik \cite{Gamarnik}. When Chv\'atal and Reed \cite{CR} formulated the threshold conjecture
in 1992, they also settled the case $K=2$. The threshold was one clause
per variable, as Goerdt \cite{Goerdt} proved independently. For larger
$K$, Friedgut \cite{Friedgut} later showed that the transition occurs
in an interval of clause densities whose width tends to zero as $N$
grows. The interval could still move as it narrowed.

Meanwhile, physicists were predicting where this interval should lie.
M\'ezard and Parisi \cite{MP1,MPzero} had used the cavity method to
propose a variational formula for the free energy of diluted spin
glasses. Their ansatz with one step of replica symmetry breaking
(1RSB) was then applied to random $K$-SAT \cite{MPZ,MZcavity}.
In this description the satisfying assignments near the threshold
fall into clusters \cite{KMRZ,MPRT,MRS}. For $3$-sat, estimating
how many clusters survive as clauses are added put the predicted
threshold at about $4.267$ clauses per variable \cite{MMZ}. Coja-Oghlan and
Panagiotou \cite{Coja2} later confirmed the large-$K$ prediction up to
an additive $o(1)$ as $K\to\infty$. Ding, Sly and Sun \cite{DSS} then
proved exact agreement for every sufficiently large fixed $K$.
For general $K$, the existence question eventually led back to
Friedgut's window. Abbe and Montanari \cite{AM} had shown that a sufficiently strong bound on its width would force its location to converge. Carenini's recent estimate \cite{Carenini} met their criterion for every fixed $K\ge3$. OpenAI \cite{OAIall} gave another proof of convergence.

The proof of the M\'ezard--Parisi formula itself drew on an idea
from the Sherrington--Kirkpatrick model. Guerra \cite{Guerra} had
used interpolation to bound the free energy of that model from above
by the Parisi variational expression. Franz and Leone \cite{FL} adapted
his interpolation to diluted models and obtained the corresponding
M\'ezard--Parisi upper bound at the replica-symmetric and one-step
levels. Once the bound had been extended to every finite depth
\cite{PT,Talagrand}, the question was whether taking the infimum
over all these trials would give the free energy exactly.
The matching lower bound for random $K$-SAT is now known for every $K\ge2$ \cite{Panchenko}, with another proof in \cite{OAImp}.

At zero temperature, Bayati, Gamarnik and Tetali \cite{BGT} had
already proved that the minimum number of violations per variable
has a limit. They also proposed using this limit to locate the
SAT threshold. In the discussion following their Conjecture~1, they
explained that exponential decay of the satisfiability probability
above the transition window would make the formula remain
unsatisfiable even after a small positive fraction of its clauses
were removed. It would follow that the SAT threshold is the largest
density at which the limiting energy vanishes. Gamarnik \cite[Section~4]{Gamarnik} later observed that their results would also show that this limit is positive above the
SAT threshold if one could prove exponential concentration for the
number of clauses at which satisfiability is first lost.\vspace{-2mm}

\paragraph{Outline of the paper.}
We state the threshold characterization in Theorem~\ref{thm:main}. Section~\ref{sec:zeroMP} starts with the positive-temperature
MP formula and then gives its zero-temperature version for the limit 
of the ground-state energy. Before proving this formula, 
we look at the one-step formula in Section~\ref{sec:one-step}
and explain why it gives the threshold of Ding, Sly and Sun for
large $K$. The derivation of the zero temperature formula is in
Section~\ref{sec:zeroMP-proof}.
Section~\ref{sec:robust} proves that zero limiting energy at one density implies satisfiability at every smaller density with high probability. The proof uses a uniform deletion estimate, which is proved in Section~\ref{sec:deletion} by adapting OpenAI's argument for $3$-SAT \cite[Lemma~4.1]{OAI3}.
Section~\ref{sec:transition} then identifies the
SAT threshold by the zero-temperature MP formula.\vspace{-2mm}

\paragraph{Acknowledgement.}
In addition to finding the arguments in the literature, ChatGPT helped me draft this survey. I rewrote many of the arguments myself in a form I found easier to follow. In doing so, I noticed that the deletion argument, together with the ground-state limit, already proves the existence of the satisfiability threshold without using Carenini's theorem as an input.

\subsection{The variational threshold}

Let $\Sigma_N=\{-1,1\}^N$.
A clause $C$ chooses $K$ indices $i_1,\ldots,i_K$ independently and
uniformly from $[N]=\{1,\ldots,N\}$, together with independent symmetric
signs $J_1,\ldots,J_K$. Repeated indices are allowed. We use the
convention that the clause is violated by an assignment vector $\sigma\in\Sigma_N$ precisely
when $\sigma_{i_j}=J_j$ for every $j$. Its violation indicator is
\[
 U_C(\sigma)=\prod_{j=1}^K\frac{1+J_j\sigma_{i_j}}2.
\]
Let $\lambda\ge0$ be the clause density and let $\pi_0$ be a rate-one
Poisson process, independent of the sequence of independent clauses
$(C_a)_{a\ge1}$. Thus $\pi_0(\lambda N)$ has mean $\lambda N$. Write
\begin{equation}\label{eq:model-quantities}
 H_N(\sigma)=\sum_{a=1}^{\pi_0(\lambda N)}U_{C_a}(\sigma),
 \qquad V_N(\lambda)=\min_{\sigma\in\Sigma_N}H_N(\sigma).
\end{equation}
The formula is satisfiable if and only if $V_N(\lambda)=0$. At
inverse temperature $\beta>0$ the MP Hamiltonian is $-\beta H_N$,
so a violation lowers the exponent in the Boltzmann weight.

Section~\ref{sec:zeroMP} defines the explicit zero-temperature
functional $\MP^{(0)}_{K,\lambda}(\mathcal T)$. Its trials, or order parameters, are finite hierarchies
\[
 \mathcal T=(r,\boldsymbol y,h),\qquad
 r\ge1,\quad 0<y_0<\cdots<y_{r-1}<\infty,\quad
 h:[0,1]^r\longrightarrow[-1,1]
\]
with $h$ measurable. The field $h$ records a signed energy gap in the M\'ezard-Parisi ansatz. Denote
the variational value and the infimum of the densities where it is negative by
\begin{equation}\label{eq:main}
 M_K^{(0)}(\lambda)=\inf_{\mathcal T}\MP^{(0)}_{K,\lambda}(\mathcal T),
 \qquad
 \alpha_K=\inf\{\lambda\ge0:M_K^{(0)}(\lambda)<0\}.
\end{equation}
The trial $h\equiv0$ has value zero, so $M_K^{(0)}\le0$.
The following theorem identifies this variationally defined $\alpha_K$ as the satisfiability threshold.

\begin{theorem}[The SAT--UNSAT transition]\label{thm:main}
For every fixed $K\ge3$, the number $\alpha_K$ in \eqref{eq:main}
is finite and strictly positive, and
\begin{equation}\label{eq:sat-transition}
 \lim_{N\to\infty}\Pp\bigl(V_N(\lambda)=0\bigr)
 =\begin{cases}
  1,&0\le\lambda<\alpha_K,\\
  0,&\lambda>\alpha_K.
 \end{cases}
\end{equation}
For each $\lambda>\alpha_K,$ there is $\varepsilon>0$ such that $\Pp(V_N(\lambda)\ge\varepsilon N)\to1$.
\end{theorem}

The link to the thermodynamic description is the ground-state identity
\begin{equation}\label{eq:energy}
 e_K(\lambda):=\lim_{N\to\infty}\frac{\E V_N(\lambda)}N
            =-M_K^{(0)}(\lambda).
\end{equation}
When $e_K(\lambda)>0$, concentration of $V_N(\lambda)$ gives a linear number of violations with high probability. When $e_K(\lambda)=0$, the question is how to pass from a vanishing fraction of violations to exact satisfiability. The deletion argument makes this possible at every strictly smaller density. At $\alpha_K$ itself, continuity gives $e_K(\alpha_K)=0$.

\section{The zero-temperature MP formula}\label{sec:zeroMP}

\paragraph{The positive-temperature MP formula.}
The model defined above is the one used in \cite{Panchenko}, with
independent variable indices and a Poisson number of clauses.
At inverse temperature $\beta>0$, a clause has Hamiltonian
\[
 \theta_\beta(\sigma_1,\ldots,\sigma_K)
   =-\beta\prod_{j=1}^K\frac{1+J_j\sigma_j}{2},
 \qquad
 H_N^{\rm model}(\sigma)=
 \sum_{a=1}^{\pi_0(\lambda N)}
       \theta_{\beta,a}(\sigma_{i_{1,a}},\ldots,\sigma_{i_{K,a}}).
\]
The partition function is
\begin{equation}\label{eq:pressure-definition}
 Z_N(\lambda,\beta)=\sum_{\sigma\in\Sigma_N}
                   e^{H_N^{\rm model}(\sigma)},
\end{equation}
and the free energy is $N^{-1}\E\log Z_N$. For spin means
$x_j\in[-1,1]$, the clause is extended by averaging its Boltzmann
factor over independent spins with those means and taking the
logarithm. The product below is the probability that all literals
are false. In the notation of \cite{Panchenko}, the extension is
\begin{equation}\label{eq:soft-clause}
 \theta_\beta(x_1,\ldots,x_K)
 =\log\Biggl(1-(1-e^{-\beta})
                   \prod_{j=1}^K\frac{1+J_jx_j}{2}\Biggr).
\end{equation}

A hierarchical trial consists of a depth $r\ge1$, parameters
$0<\zeta_0<\cdots<\zeta_{r-1}<1$, and a measurable map
$h_0:[0,1]^r\to[-1,1]$. For every input site $I$, take independent
uniform variables $Z_1^I,\ldots,Z_r^I$ and put
\begin{equation}\label{eq:MPstandard}
 m^I=h_0(Z_1^I,\ldots,Z_r^I).
\end{equation}
The map $h_0$ is the same at every site. The coordinates are
independent across sites and levels. The level index specifies when a
coordinate is averaged.

In the large-system limit a new variable is incident to $\lambda K$
clauses on average. At fixed density adding a variable increases the
total number of clauses by only $\lambda$ on average, accounting for a
correction of mean $\lambda(K-1)$. Take independent counts
$N_1\sim\Poi(\lambda K)$ and $N_2\sim\Poi(\lambda(K-1))$ for the
incident and correction clauses. Index a site by $I=(\ell,a,j)$ for
input $j$ of clause $a$ of kind $\ell$. Each gets a fresh copy of
\eqref{eq:MPstandard}. The clauses $\theta_{\beta,\ell,a}$ have
independent signs $J_{j,a}^{\ell}$, independent also of the counts
and site coordinates. With the new spin held at $\varsigma$, put
\begin{align}
 A^0(\varsigma)
 &=\sum_{a=1}^{N_1}\theta_{\beta,1,a}
        (m^{1,a,1},\ldots,m^{1,a,K-1},\varsigma),
        \qquad \varsigma=\pm1,\notag\\
 B^0
 &=\sum_{a=1}^{N_2}\theta_{\beta,2,a}
        (m^{2,a,1},\ldots,m^{2,a,K}).
                                             \label{eq:finite-exponents}
\end{align}
The two values $A^0(+1)$ and $A^0(-1)$ use the same clauses and
input means. For any terminal quantity $X$ formed from these inputs,
average its site coordinates one level at a time, starting with the last.
Write
\begin{equation}\label{eq:finite-recursion}
 X_r=X,\qquad
 X_j=\frac1{\zeta_j}\log\E_j e^{\zeta_jX_{j+1}}
       \quad(j=r-1,\ldots,0),\qquad
 \mathcal T_{\boldsymbol\zeta}(X)=X_0.
\end{equation}
Here $\E_j$ integrates the coordinates $Z_{j+1}^I$ at all input
sites, holding earlier coordinates, counts and signs fixed. The
MP functional is
\begin{equation}\label{eq:finite-functional}
 \MP_0(r,\boldsymbol\zeta,h_0)
 =\E\mathcal T_{\boldsymbol\zeta}
       \Biggl(\log\sum_{\varsigma=\pm1}e^{A^0(\varsigma)}\Biggr)
       -\E\mathcal T_{\boldsymbol\zeta}(B^0).
\end{equation}
The outer expectation averages the counts and signs. We keep the
notation $\MP_0$ of \cite{Panchenko} for the unperturbed functional.
The sum over the two spin values absorbs the $\log2$ in its original
definition \cite[Section~2.2]{Panchenko}. The recursive expression
follows from the averaging identity in \cite[Theorem~2.9]{SKbook}.

The next theorem combines the interpolation upper bound
\cite{FL,PT}, including its odd-$K$ SAT extension
\cite[Section~6.5]{Talagrand}, with the lower bound of
\cite[Theorem~2.1]{Panchenko}.

\begin{theorem}[Positive-temperature MP formula]\label{thm:finiteMP}
For every $K\ge2$, $\lambda\ge0$ and $\beta>0$, the free energy
has the limit
\begin{equation}\label{eq:MP-input}
 \begin{aligned}
 P_K(\lambda,\beta)
 &:=\lim_{N\to\infty}\frac1N\E\log Z_N(\lambda,\beta)\\
 &=\inf_{r,\boldsymbol\zeta,h_0}\MP_0(r,\boldsymbol\zeta,h_0)
   =:M_K(\lambda,\beta).
 \end{aligned}
\end{equation}
The infimum is over all finite depths, parameter vectors and
measurable field maps specified above.
\end{theorem}

\paragraph{The ground-state limit and the zero-temperature formula.}
A minimizing assignment contributes $e^{-\beta V_N}$ to the
partition function, and each of the $2^N$ assignments contributes
at most this much. Therefore,
\begin{equation}\label{eq:sandwich}
 -\frac{V_N}{N}\le\frac{1}{\beta N}\log Z_N\le \frac{\log 2}{\beta}-\frac{V_N}{N}.
\end{equation}
Taking expectations and using the existence of the free energy limit for each $\beta>0$, we see that the limsup and liminf of $\E V_N/N$ differ by at most $\log 2/\beta$. Letting
$\beta\to\infty$ proves the existence of the limit $e_K(\lambda)$,
proved for every $K\ge2$ in \cite{BGT}. Taking $N\to\infty$ in the same
bounds gives
\begin{equation}\label{eq:pressureenergy}
 - e_K(\lambda)\le \frac{P_K(\lambda,\beta)}{\beta}
                       \le\frac{\log 2}{\beta}- e_K(\lambda).
\end{equation}
In particular,
\begin{equation}\label{eq:zerolim}
 \lim_{\beta\to\infty}\frac{P_K(\lambda,\beta)}{\beta}
 =-\lim_{N\to\infty}\frac1N\E\min_{\sigma\in\Sigma_N}H_N(\sigma)
 =-e_K(\lambda).
\end{equation}
The zero temperature MP formula for $e_K(\lambda)$ is obtained by taking the $\beta\to\infty$ limit of the functional in \eqref{eq:finite-functional} after dividing the functional by \(\beta\). We describe it next.

At zero temperature, a trial consists of a depth $r\ge1$, marks $0<y_0<\cdots<y_{r-1}<\infty$, and a measurable trial map $h:[0,1]^r\to[-1,1]$. Using the same site coordinates as in
\eqref{eq:MPstandard}, put $h^I=h(Z_1^I,\ldots,Z_r^I)$.
Write $(x)_+=\max\{x,0\}$. With the counts and signs from
\eqref{eq:finite-exponents}, define
\begin{align}
 w_a^1&=\min_{1\le j<K}(J_{j,a}^1h^{1,a,j})_+,
 &w_a^2&=\min_{1\le j\le K}(J_{j,a}^2h^{2,a,j})_+,
                                      \label{eq:zerowarnings}\\
 E_1&=\min_{\varsigma=\pm1}\sum_{a=1}^{N_1}
       \one\{\varsigma=J_{K,a}^1\}w_a^1,
 &E_2&=\sum_{a=1}^{N_2}w_a^2.          \label{eq:zerocosts}
\end{align}
With the same conditional expectations $\E_j$, the recursion becomes
\begin{equation}\label{eq:energetic-recursion}
 X_r=X,\qquad
 X_j=\frac1{y_j}\log\E_j e^{y_jX_{j+1}}
       \quad(j=r-1,\ldots,0),\qquad
 \RR_{\boldsymbol y}(X)=X_0.
\end{equation}
The \emph{zero-temperature MP functional} and its variational value are
\begin{align}
 \MP^{(0)}_{K,\lambda}(r,\boldsymbol y,h)
 &=\E\RR_{\boldsymbol y}(-E_1)
          -\E\RR_{\boldsymbol y}(-E_2),\label{eq:zero-functional}\\
 M_K^{(0)}(\lambda)
 &=\inf_{r,\boldsymbol y,h}
   \MP^{(0)}_{K,\lambda}(r,\boldsymbol y,h).
                                                    \label{eq:zero-value}
\end{align}
For $h\equiv0$ both costs vanish, giving a trial of value zero at
every density.

\begin{theorem}[Zero-temperature MP formula]\label{thm:zeroMP}
For every $K\ge2$ and $\lambda\ge0$, the limit in \eqref{eq:zerolim}
is given by
\begin{equation}\label{eq:zeroMPformula}
 -e_K(\lambda)
 =M_K^{(0)}(\lambda)
 =\inf_{r,\boldsymbol y,h}\MP^{(0)}_{K,\lambda}(r,\boldsymbol y,h).
\end{equation}
The infimum is over all finite energetic trials specified above.
\end{theorem}

\begin{remark}
To interpret the quantities in \eqref{eq:zerowarnings} and \eqref{eq:zerocosts}, imagine adding a clause $C$ with $q$ literals to a finite formula. Let $V$ be the minimum number of violations before adding $C$, and let $V'$ be the
minimum afterwards. For the variable in position $j$ of $C$, let
$V_j(s)$ be the old minimum when this variable is fixed at $s$.
Thus $V=\min\{V_j(+1),V_j(-1)\}$. The signed gap
$t_j=V_j(-1)-V_j(+1)$ is positive when $+1$ is the cheaper choice.
Under our sign convention, the $j$th literal is true when its spin
equals $-J_j$. Requiring this literal to be true therefore costs
\(
 V_j(-J_j)-V=(J_jt_j)_+
\)
extra violations in the old formula.

Since any one of the $q$ literals suffices to satisfy $C$, the least
extra cost of satisfying it is $\min_{1\le j\le q}(J_jt_j)_+$.
If an old minimizing assignment already satisfies $C$, this cost is
zero. Otherwise we may also keep an old minimizing assignment and
incur one additional violation in $C$. The increase in the minimum
number of violations is therefore
\begin{equation}\label{eq:minplus-clause}
 V'-V=\min\left\{1,\min_{1\le j\le q}(J_jt_j)_+\right\}
     =\min_{1\le j\le q}(J_jh_j)_+,
\end{equation}
where $h_j=\clip(t_j)$ and $\clip(t)=\max\{-1,\min\{t,1\}\}$.

In the MP description of the thermodynamic limit, the fields $h^I$ supplied by the hierarchical trial play the role of the clipped energy gaps $h_j$ in \eqref{eq:minplus-clause}.

The \emph{site term} $E_1$ comes from adding a new spin together with its incident clauses. In the MP trial, each incident clause receives its own independent input fields. For a fixed value \(\varsigma\) of the new spin, clause \(a\) contributes \(w_a^1\) if its literal on the new spin is false, $\varsigma=J_{K,a}^1$, and zero otherwise. We add the contributions for each of the two spin values and take the smaller total, since the same spin occurs in every incident clause.

The \emph{clause correction} $E_2$
accounts for the clause density in the cavity correction term. The new spin brings in $\lambda K$
incident clauses on average. At fixed density $\lambda$, adding one
spin calls for a net increase of only $\lambda$ clauses on average.
This explains why the correction involves a Poisson number of full
clauses on the old spins with mean $\lambda(K-1)$. Each such clause
contributes $w_a^2$, and their total cost is $E_2$.

Finally, we refer to Remark \ref{rem:postempgaps} below to see how these quantities arise from the corresponding ones at positive temperature.
\qed
\end{remark}

We turn next to the one-step specialization of \eqref{eq:zeroMPformula}.
Section~\ref{sec:zeroMP-proof} then derives the clipped fields from
the rescaled log odds of the positive-temperature trial spins and
proves the zero-temperature formula.

\section{One-step RSB and the SAT threshold}\label{sec:one-step}\label{subsec:one-step}

What is usually called one-step replica symmetry breaking in the random $K$-SAT model at zero temperature is, in fact, a limiting case of the MP trials we considered above. Namely, if we take $0<y_0<y,$ take $h^I=h(Z_1^I,Z_2^I)$ and let $y_0\downarrow 0$, we get the one-step functional
\begin{equation}\label{eq:one-step-functional}
 \MP^{(0),1}_{K,\lambda}(\mathcal Q,y)
 =\frac1y\E\log\E_h e^{-yE_1}
       -\frac1y\E\log\E_h e^{-yE_2}.
\end{equation}
As \(y_0\downarrow0\), the average over the first coordinates \(Z_1^I\) becomes an ordinary expectation outside the logarithm. The conditional law \(\nu^I=\mathcal L(h^I\mid Z_1^I)\) is called the survey at site \(I\). The surveys are independent with common law \(\mathcal Q\). Given the surveys, \(\mathbb E_h\) averages the independent gaps \(h^I\) drawn from them. The outer expectation averages the surveys, clause counts and signs.

Let $M_K^{(0),1}$ be the infimum of \eqref{eq:one-step-functional}
over $\mathcal Q,y$, and set
\begin{equation}\label{eq:one-step-threshold}
 \alpha_K^{1}=\inf\{\lambda\ge0:M_K^{(0),1}(\lambda)<0\}.
\end{equation}
Restricting the full infimum to one step and then testing the zero
field gives $M_K^{(0)}\le M_K^{(0),1}\le0$ and
$\alpha_K\le\alpha_K^1$. To see
that $\alpha_K^1$ is finite, take $h\equiv1$.
By Poisson thinning the numbers of non-zero clauses in $E_1$ in \eqref{eq:zerocosts} for each of the two spin values $\varsigma=\pm 1$ are independent $X_+,X_-\sim\Poi(\mu)$,
where $\mu=\lambda K2^{-K}$. Thus $E_1=\min(X_+,X_-)$ and
$\E E_2=\lambda(K-1)2^{-K}$. The trial value is
\[
\E(-E_1+E_2) = 
 -\lambda 2^{-K}+\tfrac12\E|X_+-X_-|
 \le -\lambda 2^{-K}+\sqrt{\mu/2},
\]
which is negative for large $\lambda$.

The result of Ding, Sly and Sun \cite{DSS} implies that
\begin{equation}\label{eq:DSS-consequence}
 \alpha_K=\alpha_K^1
 \qquad\text{for every sufficiently large }K.
\end{equation}
Their threshold is the zero of the predicted \emph{complexity}, the
exponential growth rate of the number of solution clusters, computed
at a particular distributional fixed point. Their SAT lower bound
ensures zero energy below this density. At each fixed density just above the threshold, their positive-temperature one-step trial with
$\zeta=\beta^{-1/2}$ has value at most $-c\sqrt{\beta}$ for some $c>0$
and all large $\beta$ \cite[Section~4, \emph{Threshold upper bound}]{DSS}.
Applying \eqref{eq:trial-temperature} below with the surveys held fixed
during the recursion gives an energetic trial with $y=\sqrt{\beta}$
and value at most $-c/\sqrt{\beta}+C_{K,\lambda}/\beta<0$ for large $\beta$.
These two facts imply \eqref{eq:DSS-consequence} for large enough $K$.
On the other hand, the cavity prediction of \cite{MMZ} gives a threshold by the same fixed-point prescription for every $K\ge3$. A corresponding
variational question is whether $\alpha_K=\alpha_K^1$ for all such
$K\geq 3$, which remains open.

\section{Proof of the zero-temperature MP formula}\label{sec:zeroMP-proof}

The spin-mean map $h_0$ in the MP infimum may vary with $\beta$.
We compare a single clause uniformly over its input means and then
pass the estimate through the hierarchy, similarly to \cite{AC}.
Under the change of variables $y_j=\beta\zeta_j$, the condition
$\zeta_{r-1}<1$ becomes $y_{r-1}<\beta$. For $B>0$ let
$M_{K,<B}^{(0)}(\lambda)$ be the energetic infimum \eqref{eq:zero-value}
with extra constraint $y_{r-1}<B$.

\begin{proposition}[Uniform temperature comparison]\label{prop:temperature-comparison}
For every $K\ge2$, $\lambda\ge0$ and $\beta>0$,
\begin{equation}\label{eq:uniform-temperature}
 \left|\frac{M_K(\lambda,\beta)}{\beta}-M_{K,<\beta}^{(0)}(\lambda)\right|
 \le \frac{C_{K,\lambda}}{\beta},
\end{equation}
where $C_{K,\lambda}=\log 2+\lambda K\log K+\lambda(K-1)\log(K+1).$
\end{proposition}

\begin{proof}
Consider $q\ge1$ independent input spins with arbitrary means $m_j$
and prescribed literal signs $J_j$. Write $a_j=(1-J_jm_j)/2$ for
the probability that the $j$th literal is true. The averaged clause
factor in \eqref{eq:soft-clause} is
\[
 G_\beta=1-(1-e^{-\beta})\prod_{j=1}^q(1-a_j).
\]
To realize this factor as a probability, accept a clause if any
literal is true or if an independent coin of success probability
$e^{-\beta}$ succeeds. The probability of a union is bounded below
by its largest term and above by their sum. Thus, 
\[
\max\{e^{-\beta},a_1,\ldots,a_q\}\le G_\beta\le e^{-\beta}+\sum_{j=1}^q a_j\le(q+1)\max\{e^{-\beta},a_1,\ldots,a_q\}.
\]
If we let
\begin{equation}\label{eq:literal-cost}
 \ell_{\beta,J}(m)=-\frac1\beta\log\frac{1-Jm}{2},
 \qquad -\log0=+\infty,
\end{equation}
then taking logarithms gives
\begin{equation}\label{eq:clause-literal-cost}
 0\le \frac{\log G_\beta}{\beta}+
 \min\{1,\min_j\ell_{\beta,J_j}(m_j)\}
 \le\frac{\log(q+1)}{\beta}.
\end{equation}
We can now rewrite this as follows. For $-1<m<1$, the smaller value
$$
\delta_\beta(m)=\min_{J=\pm1}\ell_{\beta,J}(m)
=\min\Bigl(-\frac1\beta\log\frac{1-m}{2}, -\frac1\beta\log\frac{1+m}{2}\Bigr)
$$ 
is at most $\log 2/\beta$. Their difference is the rescaled log odds
\begin{equation}\label{eq:field-rescaling}
 t_\beta(m)=\ell_{\beta,+1}(m)-\ell_{\beta,-1}(m)
           =\frac1\beta\log\frac{1+m}{1-m}=\frac2\beta \tanh^{-1}(m).
\end{equation}
Thus $\ell_{\beta,J}(m)=\delta_\beta(m)+(Jt_\beta(m))_+$ for
either sign of $J$. Dropping the smaller cost $\delta_\beta(m_j)$ at each
site reduces $\min\{1,\min_j\ell_{\beta,J_j}(m_j)\}$ by at most $\log2/\beta$. If we define $h_\beta(m)=\clip(t_\beta(m))$ with $\clip(t)=\max\{-1,\min\{t,1\}\}$ then the remaining minimum in \eqref{eq:clause-literal-cost} is
\[
 \min\{1,\min_j(J_jt_\beta(m_j))_+\}
    =\min_j(J_jh_\beta(m_j))_+.
\]
Together with
\eqref{eq:clause-literal-cost} this gives
\begin{equation}\label{eq:clause-temperature}
 -\frac{\log2}{\beta}
 \le\frac{\log G_\beta}{\beta}+\min_j(J_jh_\beta(m_j))_+
 \le\frac{\log(q+1)}{\beta}.
\end{equation}
The estimate extends by continuity to $m_j=\pm1$, with
$t_\beta(\pm1)=\pm\infty$ and $h_\beta(\pm1)=\pm1$.

Recall that $m_j$ are generated by $h_0$ in \eqref{eq:MPstandard}, and
note that the composition $h=h_\beta\circ h_0$ gives \emph{every}
measurable gap map $h:[0,1]^r\to[-1,1]$ at each fixed $\beta>0$.
Indeed, given $h$, take
\begin{equation}\label{eq:reverse-temperature}
 h_{0}(u)=\tanh\bigl(\beta h(u)/2\bigr).
\end{equation}
Then \eqref{eq:field-rescaling} gives $t_\beta(h_{0}(u))=h(u)$ even when $h(u)=\pm1$.
Clipping does not do anything because $h$ is already in $[-1,1]$, so $h_\beta(h_{0}(u))=h(u).$ Therefore, $y_{r-1}<\beta$ is the only restriction on the energetic trials in this comparison.

Let us now approximate the two terms in \eqref{eq:finite-exponents} with the counts, signs and site coordinates fixed. A clause $\theta$ in $A^0(\varsigma)$ contributes zero whenever the new spin $\varsigma=\pm1$ satisfies its literal, $\frac{1+J\varsigma}{2}=0$. For the other spin value, the clause has \(q=K-1\) remaining inputs, so the error in approximating \(\theta_{\beta,1,a}/\beta\) by \(-w_a^1\) is at most \((\log K)/\beta\).
If we recall the definition of $E_1$ and $E_2$ in \eqref{eq:zerocosts} then
$$
\Biggl|\frac1\beta\sum_{a=1}^{N_1}\theta_{\beta,1,a}
        (m^{1,a,1},\ldots,m^{1,a,K-1},\varsigma)
  +
\sum_{a=1}^{N_1}
       \one\{\varsigma=J_{K,a}^1\}w_a^1
 \Biggr| 
\leq  \frac{N_1\log K}{\beta}.      
$$
A full clause in $B^0$ has $q=K$ inputs, so 
$$
\Biggl|\frac1\beta\sum_{a=1}^{N_2}\theta_{\beta,2,a}
        (m^{2,a,1},\ldots,m^{2,a,K})
  +
\sum_{a=1}^{N_2}w_a^2
 \Biggr| 
\leq  \frac{N_2\log(K+1)}{\beta}.
$$
Using that the logarithm of a sum of two exponentials differs from the larger exponent by at most $\log 2$, we obtain
\begin{align}
 \Biggl|\frac1\beta\log\sum_{\varsigma=\pm1}e^{A^0(\varsigma)}+E_1\Biggr|
 &\le\frac{\log 2+N_1\log K}{\beta},\quad
 \left|\frac{B^0}{\beta}+E_2\right|
 \le\frac{N_2\log(K+1)}{\beta}.
 \label{eq:local-temperature}
\end{align}

With $\boldsymbol y=\beta\boldsymbol\zeta$ we have
$\beta^{-1}\mathcal T_{\boldsymbol\zeta}(X)=\RR_{\boldsymbol y}(X/\beta)$.
Each step preserves order and commutes with adding a constant. An
error bounded by $c$ at the terminal level is therefore still bounded
by $c$ after every step, provided $c$ is fixed during the averaging.
The bounds in \eqref{eq:local-temperature} depend only on the counts.
Applying the recursion and then averaging gives
\begin{equation}\label{eq:trial-temperature}
 \left|\frac{\MP_0(r,\boldsymbol\zeta,h_0)}{\beta}
       -\MP^{(0)}_{K,\lambda}
          (r,\beta\boldsymbol\zeta,h_\beta\circ h_0)\right|
 \le\frac{C_{K,\lambda}}{\beta},
\end{equation}
since $\E N_1=\lambda K$ and $\E N_2=\lambda(K-1)$. This estimate is
uniform in the depth, marks and field map. Every energetic trial with
$y_{r-1}<\beta$ is obtained by \eqref{eq:reverse-temperature} and
$\zeta_j=y_j/\beta$, so taking infima proves \eqref{eq:uniform-temperature}.
\end{proof}

\begin{remark} \label{rem:postempgaps}
Return to the finite formula used in \eqref{eq:minplus-clause},
before the clause $C$ was added. Let $Z_j(s)$ be the sum of the
Boltzmann weights of all assignments in which the spin in position
$j$ of $C$ equals $s$. The unrestricted partition function is
$Z=Z_j(+1)+Z_j(-1)$, so the Gibbs probability of the value $s$ is $p_s=Z_j(s)/Z$. If $m_\beta$ is the Gibbs mean of this spin then $p_s=(1+s m_\beta)/2$.

Write $F=-\beta^{-1}\log Z$ for the unrestricted free energy and $F_j(s)=-\beta^{-1}\log Z_j(s)$ for the free energy with the spin fixed at $s$. The corresponding increase in the free energy is
$$
 F_j(s)-F
 =-\frac1\beta\log\frac{Z_j(s)}{Z}
  =-\frac1\beta\log p_s
 =-\frac1\beta\log\frac{1+sm_\beta}{2}.
$$
To make the $j$th literal in $C$ true, we must fix its spin
at $-J_j$. The resulting increase is exactly
\[
 F_j(-J_j)-F
 =-\frac1\beta\log\frac{1-J_jm_\beta}{2}
 =\ell_{\beta,J_j}(m_\beta).
\]
Finally, the difference between the two constrained free energies is the free energy gap
\[
 F_j(-1)-F_j(+1)
 =\frac1\beta\log\frac{Z_j(+1)}{Z_j(-1)}
 =\frac1\beta\log\frac{1+m_\beta}{1-m_\beta}
 =t_\beta(m_\beta).
\]
For this fixed formula, the same estimate as in \eqref{eq:sandwich}
gives $F_j(s)\to V_j(s)$ and $F\to V$ as $\beta\to\infty$.
Thus $t_\beta(m_\beta)$ converges to the signed energy gap
$V_j(-1)-V_j(+1)$, while $\ell_{\beta,J_j}(m_\beta)$ converges to
$V_j(-J_j)-V$, the extra cost of making the $j$th literal true.
\qed
\end{remark}

\begin{proof}[Proof of Theorem~\ref{thm:zeroMP}]
The classes $y_{r-1}<\beta$ increase to the class of all finite
energetic trials. Hence $M_{K,<\beta}^{(0)}(\lambda)\downarrow
M_K^{(0)}(\lambda)$. By Theorem~\ref{thm:finiteMP} and
\eqref{eq:uniform-temperature}, their limit is also the limit of
$P_K(\lambda,\beta)/\beta$, which equals $-e_K(\lambda)$ by
\eqref{eq:pressureenergy}.
\end{proof}

\begin{corollary}[Finite energetic caps]\label{cor:energetic-cap}
For every $B>0$,
\begin{equation}\label{eq:energetic-cap}
 0\le M_{K,<B}^{(0)}(\lambda)-M_K^{(0)}(\lambda)
       \le\frac{C_{K,\lambda}+\log2}{B}.
\end{equation}
\end{corollary}
\begin{proof}
The lower bound follows by inclusion of the trial classes. For the
upper bound, take $\beta=B$ in \eqref{eq:uniform-temperature} and
use $P_K(\lambda,B)/B\le M_K^{(0)}(\lambda)+\log2/B$ from
\eqref{eq:pressureenergy}.
\end{proof}

\section{From zero energy to satisfiability}\label{sec:robust}

We compare the formula at two densities $a<b$ by following the same sequence of clauses. Suppose the minimum number of violations at density $b$ is $o(N)$ in expectation. Could the formula already have been unsatisfiable at density $a$?

\begin{proposition}[Satisfiability at a smaller density]\label{prop:robust}
For every fixed \(K\ge3\) and \(0\le a<b\),
\begin{equation}\label{eq:goal}
 e_K(b)=0 \quad\Longrightarrow\quad
 \Pp\bigl(V_N(a)=0\bigr)\longrightarrow1.
\end{equation}
\end{proposition}

We will treat the clauses as arriving in a rate-one Poisson process, so their number at time \(\lambda N\) has law
\(\Poi(\lambda N)\). Let \(T_0\) be the first time at which no assignment satisfies all the clauses that have arrived. Thus the formula at time \(\lambda N\) is satisfiable exactly when \(T_0>\lambda N\).

\paragraph{Outline of the proof.}
The question about a linear number of violations leads us to compare two times. 
Consider a realization in which the formula is unsatisfiable at time \(aN\), yet some assignment violates at most \(r\) clauses at time \(bN\). Deleting one variable from each violated clause would leave that assignment satisfying every remaining clause at \(bN\). In this realization, therefore, the freedom to delete up to \(r\) variables must delay the first failure by more than \((b-a)N\).
As the clauses arrive, \(F_0\) will record the time when the
formula first becomes unsatisfiable, and \(F_r\) the time when even
removing up to \(r\) variables of our choice together with their
incident clauses no longer restores satisfiability. We stop both clocks at a fixed time $BN>bN$. If the original formula was already unsatisfiable at $aN$, the deletions would have to buy us at least $(b-a)N$ extra time.

We will telescope the delay as
\(
 F_r-F_0=\sum_{j=1}^r(F_j-F_{j-1})
\)
and, to control each increment, we will use related gains
\(F_{j-1}^{+v}-F_{j-1}\) as \(v\) ranges over the variables.
In \(F_{j-1}^{+v}\) we remove one prescribed variable \(v\) and are
still free to choose the other \(j-1\) deletions. The key estimate in the proof will
show that 
\[
 \Pp(F_{j-1}^{+v}-F_{j-1}>u)\le C(1+u)^{-2(K-1)/K}
\]
uniformly in \(N,j,v\). Thus the gain has a bounded \(p\)th moment for
\(1<p<2(K-1)/K\). Summing the increments will then give
\[
 \frac{\E(F_r-F_0)}{N}\le C\left(\frac rN\right)^{1-1/p}.
\]
For \(r=\lfloor\varepsilon N\rfloor\), Markov's inequality bounds the probability of a delay of \((b-a)N\) by a quantity that tends to zero with \(\varepsilon\). Now suppose \(e_K(b)=0\). For each fixed \(\varepsilon>0\), with high probability there is an assignment at \(bN\) violating at most \(r=\lfloor\varepsilon N\rfloor\) clauses. On this event, if the formula was already unsatisfiable at \(aN\), the deletions must have bought us more than \((b-a)N\) extra time. By Markov’s inequality, such probability is small and, as result, the formula is satisfiable at \(aN\) with high probability.

Most of the work goes into the tail estimate, which we prove in
Section~\ref{sec:deletion}. Suppose \(F_r^{+v}>F_r\); this means that at time
\(F_r\), removing \(v\) and some set \(D\) of at most \(r\) other
variables still leaves a satisfiable formula. We call each such
choice of \(D\), together with a satisfying assignment \(\sigma\)
of the remaining variables, a \emph{candidate}. The family \(U\)
of candidates is therefore nonempty. To bound \(F_r^{+v}-F_r\), we ask how long the last of these candidates can survive as new clauses arrive. To answer this, we will compare the probabilities
\[
 p_j(x)=\Pp\bigl(\text{a list of }\Poi(x)\text{ independent }
                       j\text{-clauses eliminates }U\bigr)
\]
for $j=K-1$ and $j=K.$
Here the clauses are sampled afresh on the variables other than \(v\).
They eliminate \(U\) if, for every candidate \((D,\sigma)\), at
least one clause avoids \(D\) and is violated by \(\sigma\).

The shorter $(K-1)$-clauses enter when we put \(v\) back. In a clause
containing \(v\) once, one value of \(v\) satisfies the clause and
the other leaves a \((K-1)\)-clause. At time \(F_r\), putting \(v\)
back with either value has already failed for every candidate.
Lemma~\ref{lem:failure-candidates} uses this fact to show that
\(F_r^{+v}>F_r\) rarely occurs with \(p_{K-1}(L_0)\) small, where
\(L_0\) is some fixed mean. Looking forward from
\(F_r\), Lemma~\ref{lem:candidate-survival} shows that the larger
\(p_{K-1}(L_0)\) is, the sooner fresh \(K\)-clauses are likely to
eliminate \(U\). Together these two estimates give the tail bound
for \(F_r^{+v}-F_r\).

\subsection{Reduction to a uniform deletion estimate}\label{sec:reduction}

Suppose an assignment violates at most \(r\) clauses. Choose one
variable from each violated clause and delete the chosen variables
together with all clauses in which they occur. The restriction of
the assignment satisfies every clause left over. Thus, to prove that
every assignment violates more than \(r\) clauses, it suffices to
show that even the best choice of \(r\) variable deletions leaves
an unsatisfiable formula.

Write \(\Phi_t\) for the original formula at time \(t\). For
\(D\subseteq[N]\), let \(T(D)\) be the first time the clauses of
\(\Phi_t\) disjoint from \(D\) have no satisfying assignment on
\([N]\setminus D\), with \(T(D)=\infty\) if this never happens. 
Notice that $T(D)\leq T(D')$ if $D\subseteq D'.$
Put
\[
 T_r=\max_{|D|\le r}T(D).
\]
The maximum is over all choices of the deleted variables and can
depend on the formula. Since clauses only arrive, \(T_r>t\) means
that some choice still works at time \(t\). For \(r=0\) the only
choice is \(D=\varnothing\), giving \(T(\varnothing)=T_0\). The
deterministic observation above gives
\begin{equation}\label{eq:repair-implication}
 \{V_N(\lambda)\le r\}\ \subseteq\ \{T_r>\lambda N\}.
\end{equation}
To make this comparison in expectation, fix an observation horizon
\(BN\) and write \(F_r=\min(T_r,BN)\). We record \(BN\) for a
formula that survives to the end of observation. This includes, for
example, the case \(r=N\), in which \(T_r=\infty\).
For our two densities we choose \(B>b\). If \(V_N(a)>0\), then \(F_0\le aN\).
If \(V_N(b)\le r\), then \eqref{eq:repair-implication} gives \(T_r>bN\), so \(F_r>bN\) as well. Thus
\begin{equation}\label{eq:repair-delay}
 \{V_N(a)>0,\ V_N(b)\le r\}\ \subseteq\
 \{F_r-F_0>(b-a)N\}.
\end{equation}
All deletion estimates below hold for each fixed \(B>0\).

To bound this delay, we increase the deletion
budget one variable at a time and write
\[
 F_r-F_0=\sum_{j=1}^r(F_j-F_{j-1}).
\]
Consider one increment, and choose a set \(S\) of \(j\) variables
whose deletion attains 
\[
F_j=\min(T_j,BN) =\min\left(\max_{|S|\le j}T(S),BN\right).
\]
We may enlarge a smaller maximizing
set to size \(j\), since deleting more variables only removes clauses.
If we prescribe any one member \(v\in S\) as a deletion, we can still
attain \(F_j\) by choosing \(S\setminus\{v\}\) as the other \(j-1\)
deletions. This suggests measuring how much one specified \(v\) adds
when we may still choose the other deletions optimally.

For a specified \(v\) and a budget \(r\) of other deletions, the
failure time and its value observed up to \(BN\) are
\[
 T_r^{+v}=\max_{\substack{D\subseteq[N]\setminus\{v\}\\ |D|\le r}}
               T(D\cup\{v\}),\qquad
 F_r^{+v}=\min(T_r^{+v},BN).
\]
Notice that $F_r\le F_r^{+v}$, and that $F_r^{+v}\le F_{r+1}$ for $r<N$.
For the first inequality, if a deletion set \(D\) in \(F_r\) contains $v$, we can still delete it by choosing $D\setminus\{v\}$ as the other deletions. Otherwise we can choose $D$ as the other deletions, and $T(D)\leq T(D\cup\{v\})$. For the second inequality, for every deletion set $D$ in $F_r^{+v}$, the union $D\cup\{v\}$ is a deletion set in $F_{r+1}$.

Moreover, if a set $S$ of cardinality $j$ attains the maximum in $F_j$, then
\begin{equation}\label{eq:optimal-members}
 F_{j-1}^{+v}-F_{j-1}=F_j-F_{j-1}\qquad \mbox{ for all } v\in S.
\end{equation}
We will bound the moment of the gain for each fixed label $v$ and sum
over all \(N\) labels. The \(j\) equal terms in \eqref{eq:optimal-members} will then bound the \(j\)th increment $F_j-F_{j-1}$.

\begin{proposition}[Uniform deletion estimate]\label{prop:moment}
For fixed \(K,B\) there is a constant \(C\) such that for all
\(N\ge1\), \(0\le r\le N\), \(v\in[N]\) and \(u\ge0\),
\begin{equation}\label{eq:deletion-tail}
 \Pp(F_r^{+v}-F_r>u)\le C(1+u)^{-2(K-1)/K}.
\end{equation}
In particular, for every
\begin{equation}\label{eq:p-range}
 1<p<2(K-1)/K
\end{equation}
there is a constant \(C_{K,B,p}\) such that
\begin{equation}\label{eq:moment}
 \E(F_r^{+v}-F_r)^p\le C_{K,B,p}.
\end{equation}
\end{proposition}

Section~\ref{sec:deletion} below proves this crucial estimate. 
One may take \(p=6/5\) for every \(K\ge3\).
The exponent $2(K-1)/K$ already appears in the proof of
Carenini's Corollary~5.3 \cite{Carenini}. For one prescribed deletion,
she combined her short-clause estimate in Lemma~2.1 with the square
in her identity~(3.4), which comes from testing both values of the
deleted variable. Since $2(K-1)/K>1$ for $K\ge3$, she obtained a
uniform bound on the mean delay. For $K=3$, OpenAI bounded the
moment of order $6/5$ uniformly in the number of other deletions allowed
\cite[Lemma~4.1]{OAI3}. Their Proposition~4.6 then controlled the
delay from an optimal set of deletions and gave a linear lower bound
on violations. We follow that argument below.

We now sum the increments to bound the mean delay.

\begin{lemma}[Delay from an optimal deletion set]\label{lem:manydeletions}
For \(p\) in \eqref{eq:p-range} and \(0\le r\le N\),
\begin{equation}\label{eq:manydeletions}
 \frac{\E F_r-\E F_0}{N}\le C_*(r/N)^{1-1/p}.
\end{equation}
\end{lemma}
\begin{proof}
For \(1\le j\le r\), sum the \(p\)th powers of the gains from
prescribing \(v\) over all \(N\) labels. By \eqref{eq:optimal-members},
at least \(j\) summands equal \((F_j-F_{j-1})^p\), so
\[
 j(F_j-F_{j-1})^p
       \le\sum_{v=1}^N(F_{j-1}^{+v}-F_{j-1})^p .
\]
Each summand now concerns a specified variable and has expectation at
most \(C_{K,B,p}\) by \eqref{eq:moment}. Taking expectations and
then \(p\)th roots gives
\(\E(F_j-F_{j-1})\le(C_{K,B,p}N/j)^{1/p}\) by Jensen's inequality.
Summing these increments and using
\(\sum_{j=1}^r j^{-1/p}\le \frac{p}{p-1} r^{1-1/p}\) proves the bound.
\end{proof}

\begin{proof}[Proof of Proposition~\ref{prop:robust}]
Fix \(B>b\) and \(p\) in \eqref{eq:p-range}. For every integer \(0\le r\le N\), \eqref{eq:repair-delay}, Markov's inequality and Lemma~\ref{lem:manydeletions} give
\begin{equation}\label{eq:two-densities}
 \Pp\bigl(V_N(a)>0\bigr)
 \le \Pp\bigl(V_N(b)>r\bigr)
       +\frac{C_*}{b-a}\left(\frac rN\right)^{1-1/p}.
\end{equation}
Suppose \(e_K(b)=0\) and let \(r=\lfloor\varepsilon N\rfloor\), with \(0<\varepsilon<1\) fixed.
Since \(\E V_N(b)=o(N)\), the first probability on the right tends to zero.
Taking the limsup in \(N\) leaves at most \(C_*\varepsilon^{1-1/p}/(b-a)\).
We can now let \(\varepsilon\downarrow0\), since \(p>1\).
\end{proof}

\section{Proof of the uniform deletion estimate}\label{sec:deletion}

Fix the variable \(v\) and the deletion budget \(r\). To prove
Proposition~\ref{prop:moment}, we stop at \(F_r=\min(T_r,BN)\). If
\(F_r=BN\), the observed delay is zero. Otherwise the original formula
has just failed even with \(r\) deletions. Deleting \(v\) in addition
may still leave satisfying assignments. We ask how long at least one
of them can continue to satisfy the new clauses that avoid \(v\).
We take \(N\ge2\), since the delay is at most \(B\) for \(N=1\).

Let \(S=[N]\setminus\{v\}\). The clauses avoiding \(v\) form the
\emph{base}. By Poisson thinning they arrive at rate
\(\rho_N=(1-1/N)^K\ge2^{-K}\), independently of the clauses
containing \(v\). A base clause has \(K\) independent uniform indices
in \(S\) and independent fair signs.

\paragraph{Candidates and their extensions.}
After deleting \(v\), we may still delete a set \(D\subseteq S\) of
at most \(r\) variables and assign values \(\sigma\in\{-1,1\}^{S\setminus D}\)
to the rest. A clause on \(S\) \emph{rejects} the pair \((D,\sigma)\)
if it avoids \(D\) and all its literals are false under
\(\sigma\). A clause meeting \(D\) is discarded. A collection of clauses
\emph{eliminates} a family of pairs if it rejects each pair at least once. Let \(U_t\) be the
set of pairs rejected by no base clause up to and including time \(t\).
These pairs are called the \emph{candidates} at \(t\). The sets of candidates \(U_t\) decrease, and
\[
 U_t\ne\varnothing\quad\Longleftrightarrow\quad T_r^{+v}>t.
\]
There are two ways to put \(v\) back while keeping \(D\) deleted.
For \(s\in\{-1,1\}\) the extension \(\sigma^{(s)}\) is
the assignment on \([N]\setminus D\) given by
\[
 \sigma^{(s)}_i=\sigma_i\quad(i\in S\setminus D),\qquad
 \sigma^{(s)}_v=s.
\]
If one extension of a pair in \(U_t\) satisfies the original clauses
through \(t\) that are disjoint from \(D\), then \(T_r>t\). Indeed,
\(D\) itself is an allowed deletion set in the original formula.
Write 
\begin{equation}
A=\{F_r<BN,\ U_{F_r}\ne\varnothing\}=\{F_r^{+v}>F_r\}.
\label{eq:eventA}
\end{equation}
This is exactly the event that deleting \(v\) buys extra time. The following property holds.

\medskip
\noindent
\textbf{Property E.} \emph{On the event $A$, for every pair $(D,\sigma)$ in \(U_{F_r}\), $\sigma$ still satisfies the base clauses disjoint from $D$, but each of its two extensions \(\sigma^{(s)}\) violates at least one original clause present at \(F_r\) and disjoint from $D$.}

\medskip
Which clause could an extension violate? First of all, such a clause must contain \(v\). If \(v\) occurs exactly once in this clause, then $v=s$ makes the literal on $v$ false and the \((K-1)\)-clause on \(S=[N]\setminus\{v\}\) obtained by erasing this literal is violated by \(\sigma\).
By time \(BN\) the expected number of clauses containing \(v\) is at
most \(L_0=KB\), since each clause contains \(v\) with probability at
most \(K/N\). For a set \(U\) of candidate pairs, let 
\begin{equation}
q(U)=\Pp\Bigl(\text{a list of }\Poi(L_0)\text{ independent } (K-1)\text{-clauses on }S\text{ eliminates } U\Bigr),
\label{eq:qU}
\end{equation}
with \(q(\varnothing)=1\). Property E gives a bound on how small \(q(U_{F_r})\) can be at failure.

\begin{lemma}\label{lem:failure-candidates}
For \(z\ge0\), with \(c_0=B\binom K2\),
\begin{equation}\label{eq:incident}
 \Pp(A,\ q(U_{F_r})\le z)\le (z+c_0/N)^2 .
\end{equation}
\end{lemma}

\begin{proof}
Let \(Q=q(U_{F_r})\).
Fix the entire base (clauses avoiding $v$) through time \(BN\). The time \(F_r\) also depends on the clauses containing \(v\), but the list of nested sets \(U_t\), \(0\le t\le BN\), is now fixed. Since the sets $U_t$ decrease with time, the probability $q(U_t)$ can only increase.
Among nonempty sets with \(q(U_t)\le z\), take the smallest one, \(W\). If there are none, the event on the
left of \eqref{eq:incident} is empty. 

On the event \(A\cap\{Q\le z\}\),
\(Q=q(U_{F_r})\leq z\), so \(W\subseteq U_{F_r}\). The advantage of \(W\) is that we chose it from the base alone, while $F_r$ in $U_{F_r}$ depends on all original clauses. Because of this, we will be able to test the two signs $s=\pm 1$ against the same fixed family and multiply the probabilities below.
 
Fix \((D,\sigma)\in W\). By Property E above, each extension \(\sigma^{(s)}\) must violate an original clause containing \(v\) and disjoint from \(D\) by time \(F_r\). In particular, for each sign $s=\pm 1$, the clauses containing \(v\) through \(BN\) must rule out every pair in \(W\).

For a given \(s\), a clause containing \(v\) can be violated by
\(\sigma^{(s)}\) only if every occurrence of \(v\) in that clause
has a literal false at \(s\). The clauses meeting this condition for
\(s=+1\) and \(s=-1\) form disjoint Poisson processes $\Pi_{+1}$ and $\Pi_{-1}$, independent of one another and of the base. Write \(E_s(W)\) for the event that, for every \((D,\sigma)\in W\), some clause of \(\Pi_s\) through \(BN\) is disjoint from \(D\) and violated by \(\sigma^{(s)}\). Since \(W\) is fixed once we condition on the base,
\begin{align*}
 \Pp(A,\ Q\le z \mid \text{base})
 &\leq \Pp(E_{+1}(W)\cap E_{-1}(W)\mid\text{base})
=\prod_{s=\pm 1} \Pp(E_s(W)\mid\text{base}).
\end{align*}
In the process $\Pi_s$, the clauses through \(BN\) containing \(v\)
exactly once give independent \((K-1)\)-clauses on \(S\) after the
false literal is erased. Their number has Poisson mean
\(\tfrac12 BK(1-1/N)^{K-1}\le L_0\). Adding independent clauses to
reach mean \(L_0\) shows that the chance they eliminate \(W\) is at
most \(q(W)\). The chance that this same process contains any clause
with \(v\) in two or more positions is at most
\(BN\binom K2/N^2=c_0/N\), by a union bound over pairs of positions.
Thus the conditional probability of ruling out every extension with
value \(s\) is at most \(q(W)+c_0/N\). Multiplying the bounds for
the two independent processes and using \(q(W)\le z\) proves the lemma.
\end{proof}

We now return to the time \(F_r\). The candidates in \(U_{F_r}\)
have survived all base clauses up to that time. How long can any
one of them continue to survive the new base clauses? Write
\(\mathcal F_t\) for the history of the full clause process through
time \(t\), and let \(\Delta=F_r^{+v}-F_r\). Given the history $\mathcal F_{F_r}$ up to the stopping time $F_r,$ suppose the event $A\in \mathcal F_{F_r}$ in \eqref{eq:eventA} holds. In this case, what does it mean that $\Delta>u$ for \(u\ge0\)?

First of all, we must have $F_r+u<BN$, because both $F_r^{+v},F_r\leq BN$. Moreover, at least one candidate in $U_{F_r}$ must survive all the base clauses arriving in $(F_r,F_r+u]$. By the strong Markov property, conditional on \(\mathcal F_{F_r}\), these $K$-clauses on \(S=[N]\setminus\{v\}\) have their original law and their number is \(\Poi(\rho_N u)\), where \(\rho_N=(1-1/N)^K\). In the notation of the next subsection, the probability that these clauses eliminate $U_{F_r}$ is denoted by $p_K(\rho_Nu)$.

Recall the notation \(Q=q(U_{F_r})\) in \eqref{eq:qU}, which is the probability that a \(\Poi(L_0)\) number of independent \((K-1)\)-clauses on \(S=[N]\setminus\{v\}\) eliminate the same set \(U_{F_r}\). In the notation of the next subsection, this probability is \(p_{K-1}(L_0)\).

Lemma~\ref{lem:candidate-survival} below relates these two probabilities as follows.

\begin{lemma}\label{lem:future-base}
Let \(\kappa=K/(K-1)\) and put \(Q=q(U_{F_r})\) as in \eqref{eq:qU}. For every \(u\ge0\),
\begin{equation}\label{eq:future-base}
 \Pp(\Delta>u\mid\mathcal F_{F_r})
 \le 2\,\one_A\exp\!\left[-\rho_N u
                         \left(\frac{Q}{2L_0}\right)^\kappa\right].
\end{equation}
\end{lemma}

On \(A\), the conditional probability of \(\Delta>u\) is at most
\(1-p_K(\rho_Nu)\). For \(u>0\), take \(U=U_{F_r}\) on the \(n=N-1\) variables of \(S\), \(x=L_0\) and \(y=\rho_Nu\) in Lemma~\ref{lem:candidate-survival}. This gives \eqref{eq:future-base}. Off \(A\) the delay is zero, and the bound at \(u=0\) follows from \(\{\Delta>0\}=A\). Using this together with Lemma~\ref{lem:failure-candidates} above, we can now prove our key deletion estimate in Proposition~\ref{prop:moment}.

\begin{proof}[Proof of Proposition~\ref{prop:moment}]
Taking expectations in
\eqref{eq:future-base} and using \(\rho_N\ge2^{-K}\) gives
\begin{equation}\label{eq:delay-inspection}
 \Pp(\Delta>u)\le 2\,\E\big[\one_A e^{-cu Q^\kappa}\big],
 \qquad c=2^{-K}(2L_0)^{-\kappa}.
\end{equation}
For \(u>0\), write the exponential as the integral of \(e^{-s}\)
from \(cuQ^\kappa\) to infinity. Then
\begin{align*}
 \E\big[\one_A e^{-cuQ^\kappa}\big]
 &=\int_0^\infty e^{-s}\,
     \Pp\big(A,\ Q\le (s/(cu))^{1/\kappa}\big)\,ds\\
 &\le \int_0^\infty e^{-s}
       \big((s/(cu))^{1/\kappa}+c_0/N\big)^2\,ds
 \le C_{K,B}\big(u^{-2/\kappa}+N^{-2}\big).
\end{align*}
The first inequality is \eqref{eq:incident}. For the second, use
\((a+b)^2\le2a^2+2b^2\) and
\(\int_0^\infty e^{-s}s^{2/\kappa}\,ds<\infty\).
For \(1\le u\le BN\), we have
\(N^{-2}\le B^2u^{-2}\le B^2u^{-2/\kappa}\). The delay is at most
\(BN\), and \(u<1\) is covered by increasing the constant. This
proves \eqref{eq:deletion-tail} with exponent \(2/\kappa=2(K-1)/K\).
Finally, integrating the tail against \(p u^{p-1}\,du\) proves
\eqref{eq:moment} for every \(p\) in \eqref{eq:p-range}.
\end{proof}

\subsection{Comparing short and long clauses}

Fix \(n\ge1\). A partial assignment is a pair \((D,\sigma)\), where \(D\subseteq[n]\) and \(\sigma\) assigns either true or false to each variable outside \(D\). A clause \emph{rejects} the pair if all its variables
lie outside \(D\) and all its literals are false under \(\sigma\).
A collection of clauses \emph{eliminates} a family \(U\) of pairs if
every member of \(U\) is rejected by at least one clause.
We sample a \(j\)-clause by drawing its \(j\) literals independently
and uniformly from \(x_1,\neg x_1,\ldots,x_n,\neg x_n\).

Fix any nonempty family \(U\) of partial assignments and, given $x\ge0$, let
$$
p_j(x)=\Pp\Bigl(\text{a list of }\Poi(x)\text{ independent }j\text{-clauses eliminates } U\Bigr).
$$
The comparison in the next lemma goes back to Friedgut's
half-cube lemma \cite[Lemma~5.7]{Friedgut}. Carenini developed the
quantitative comparison by replacing one short clause at a time
\cite[Theorem~4.3 and Lemma~6.1]{CareniniWindow}. In her threshold
paper she used Poisson interpolation to prove the counterpart of
\eqref{eq:replacement} for clauses with distinct indices
\cite[equation~(2.9)]{Carenini}. She also used the fact that survival
through two independent groups of clauses requires survival through
each group. Here we apply both ideas to families of partial assignments.
Sampling the literal positions independently gives an exact power for
the probability that one clause rejects every surviving pair, as in
\cite[proof of Lemma~4.1]{OAI3}. For the corresponding comparison
with distinct indices and full assignments, see also
\cite[Lemma~2.2 and Proposition~2.3]{OAIvariance}.

\begin{lemma}\label{lem:candidate-survival}
For any \(K\ge2\), \(x,y>0\), 
\begin{equation}\label{eq:survival}
1-p_K(y)\le 2\exp\!\left[-y\left(\frac{p_{K-1}(x)}{2x}\right)^{K/(K-1)}\right].
\end{equation}
\end{lemma}
\begin{proof}
\emph{Step 1.}
Given \(y\ge x\), we first show that
\begin{equation}\label{eq:replacement}
p_K(y)\ge p_{K-1}(x)-x(x/y)^{K-1}.
\end{equation}
Generate a list of \(M\sim\Poi(y)\) independent \(K\)-clauses $\{c_1,\ldots,c_M\}$.
At each entry, independently with probability \(a=x/y\) keep
the first \(K-1\) literals as a short clause. Otherwise leave the
entry blank. Denote these new clauses by $\{s_1,\ldots,s_M\}.$ 
Some of them could be blank and, in fact, by thinning, the resulting list contains a \(\Poi(x)\) number of independent non-blank \((K-1)\)-clauses.

Let \(\mathcal L_i = \{c_1,\ldots,c_i,s_{i+1},\ldots, s_M\}\). Thus \(\mathcal L_0\) is the thinned list and \(\mathcal L_M\) is the full list. If we condition on \(M\) first and then take expectation $\E_M$ with respect to $M$, we can write
\begin{align*}
p_K(y)- p_{K-1}(x)
&=
\E_M\Bigl[ \Pp(\mathcal L_M \text{ eliminates }U\mid M)
- \Pp(\mathcal L_0 \text{ eliminates }U\mid M)\Bigr]
\\
&=
\E_M \sum_{i=1}^M \Bigl[ \Pp(\mathcal L_i \text{ eliminates }U\mid M)
- \Pp(\mathcal L_{i-1} \text{ eliminates }U\mid M)\Bigr].
\end{align*} 
To compare the two probabilities in the last term,  let $\mathcal L_i^{-} = \{c_1,\ldots,c_{i-1},s_{i+1},\ldots, s_M\}$ be the list of their common clauses. Let us condition on $\mathcal L_i^{-}$ first, and let $W\subseteq U$ be the set of all partial assignments that survive $ \mathcal L_i^{-}.$ If \(W\) is empty, both lists already eliminate \(U\) and the conditional probabilities cancel. Otherwise, if \(W\) is non-empty, then $\mathcal L_i$ eliminates $U$ exactly when $c_i$ eliminates $W$, and $\mathcal L_{i-1}$ eliminates $U$ exactly when $s_i$ eliminates $W$.

Write \(p_j(W)\) for the probability that a single random
\(j\)-clause rejects every pair \((D,\sigma)\) in \(W\). To compute \(p_j(W)\), let \(d\) count the variables that every pair in \(W\) leaves undeleted and assigns the same value. A clause rejects all of \(W\) precisely when, at each of its \(j\) positions, the variable is among these \(d\) (with probability $d/n$) and the literal is false at the agreed value of that variable (with probability $1/2$). Therefore,
$p_j(W)=\left({d}/{2n}\right)^j$
and
$$
\Pp( c_i \text{ eliminates }W\mid \mathcal L_i^{-}, M)
- \Pp( s_i \text{ eliminates }W\mid \mathcal L_i^{-}, M)
=
\left(\frac{d}{2n}\right)^K- a \left(\frac{d}{2n}\right)^{K-1},
$$
because the clause $s_i$ is not blank with probability $a=x/y$.
If \(d/(2n)\ge a\), this change is nonnegative. If \(d /(2n)<a\), the decrease is at most \(a (d /(2n))^{K-1}\le a^K\).  Averaging over the
clauses $\mathcal L_i^{-}$ therefore gives
\[
\Pp(\mathcal L_i \text{ eliminates }U\mid M)
- \Pp(\mathcal L_{i-1} \text{ eliminates }U\mid M)
\ge -a^K.
\]
Plugging this into the above sum and taking expectations in \(M\) proves \eqref{eq:replacement}.

\emph{Step 2.}
If \(p_{K-1}(x)=0\), then \eqref{eq:survival} follows from \(1-p_K(y)\le1\). We may therefore assume \(p_{K-1}(x)>0\). Let us take $y=y_0$ in \eqref{eq:replacement} such that
\begin{equation}\label{eq:ell}
\frac{p_{K-1}(x)}{2} = x(x/y_0)^{K-1}
\Longleftrightarrow
y_0= x(2x/p_{K-1}(x))^{1/(K-1)}.
\end{equation}
We know that $p_{K-1}(x)\leq 1-e^{-x}\leq x$, because eliminating the nonempty family requires at least one clause and, in particular, $y_0\geq x.$ As a result, \eqref{eq:replacement} holds and implies that
$$
p_K(y_0)\ge p_{K-1}(x)-x(x/y_0)^{K-1} = \frac{p_{K-1}(x)}{2}.
$$
In other words, a list of \(\Poi(y_0)\) independent \(K\)-clauses has probability at least \({p_{K-1}(x)}/{2}\) of eliminating the original family \(U\).

Two independent groups with Poisson means \(s\) and \(t\) together
have the law of a single group with Poisson mean \(s+t\). If a pair
survives their union, that same pair survives each group. Hence,
by independence,
\(
 1-p_K(s+t)\le (1-p_K(s))(1-p_K(t)).
\)
Given $\Poi(y)$ clauses, if we take \(\lfloor y/y_0\rfloor\) disjoint groups of mean \(y_0\) and discard the remaining clauses,
\begin{align*}
 1-p_K(y)&\le (1-p_K(y_0))^{\lfloor y/y_0\rfloor}
\le (1-p_{K-1}(x)/2)^{\lfloor y/y_0\rfloor}
\le e^{p_{K-1}(x)/2}\exp\!\left(-\frac{p_{K-1}(x)}{2y_0}y\right).
\end{align*}
Now \(e^{p_{K-1}(x)/2}<2\) and \(p_{K-1}(x)/(2y_0)=(p_{K-1}(x)/(2x))^{K/(K-1)}\), as claimed.
\end{proof}

\section{Identifying the variational threshold}\label{sec:transition}

Fix \(K\ge3\). The zero-temperature MP formula of Theorem~\ref{thm:zeroMP}
gives the limiting mean number of violations per variable,
\[
 e_K(\lambda):=\lim_{N\to\infty}\frac{\E V_N(\lambda)}{N}
 =-M_K^{(0)}(\lambda),\qquad
 M_K^{(0)}(\lambda)=\inf_{\mathcal T}\MP^{(0)}_{K,\lambda}(\mathcal T).
\]
Here \(\mathcal T\) ranges over finite hierarchical trials.
By \eqref{eq:main}, \(\alpha_K\) is the infimum of the densities at which \(e_K\) is positive.

\begin{proof}[Proof of Theorem~\ref{thm:main}]
We couple two densities by using the same clause process. Each extra clause
increases the optimum by at most one. Taking expectations and limits
gives
\begin{equation}\label{eq:energy-Lipschitz}
 0\le e_K(\lambda')-e_K(\lambda)\le\lambda'-\lambda
 \qquad(0\le\lambda\le\lambda').
\end{equation}
Thus \(e_K\) is nondecreasing and continuous, with \(e_K(0)=0\). The trial \(h\equiv1\) in Section~\ref{sec:one-step} has negative value for large \(\lambda\), so \(\alpha_K<\infty\).

At every fixed density the optimum concentrates exponentially about
its mean on the scale $N$. Fix $\lambda>0$, let $M=\pi_0(\lambda N)$
and $m=\lfloor\lambda N\rfloor$, and write $V_{N,m}$ for the minimum
number of violations among the first $m$ clauses. Replacing one clause
changes $V_{N,m}$ by at most one, so the bounded differences inequality
gives, for $m\ge1$ and $t>0$,
\[
 \Pp\bigl(|V_{N,m}-\E V_{N,m}|\ge t\bigr)
 \le 2\exp(-2t^2/m).
\]
Using the same clause sequence gives
$|V_N(\lambda)-V_{N,m}|\le |M-m|$, so their means differ by at most
$\E|M-m|\le\sqrt{\lambda N}+1$. Combining the preceding inequality
with the Poisson tail bound for $M$, we obtain, for every $\delta>0$
and all sufficiently large $N$,
\[
 \Pp\bigl(|V_N(\lambda)-\E V_N(\lambda)|\ge\delta N\bigr)
 \le e^{-c_{\lambda,\delta}N}
\]
for some $c_{\lambda,\delta}>0$.
Consequently \(V_N(\lambda)/N\to e_K(\lambda)\) in probability.

Next, let us note that $\alpha_K\ge1$. Fix $0<\lambda<1$. The expected
number of clauses whose $K$ indices all coincide is
$\lambda N^{2-K}=o(1)$, so with high probability there are no such
clauses. In every other clause, retain the first literal and the first
subsequent literal on a different variable. By symmetry and Poisson
thinning, the retained pairs form a Poisson sample of uniform proper
$2$-clauses with limiting density $\lambda$. The subcritical $2$-SAT
theorem \cite{BBCKW} therefore gives satisfiability of the shortened
formula with high probability. On the event that no clause was
discarded, any assignment satisfying the shortened formula also
satisfies the original one. Thus $V_N(\lambda)=0$ with high
probability, and the convergence above gives $e_K(\lambda)=0$,
as claimed.

 We can now read off the transition. By \eqref{eq:energy-Lipschitz} and the definition of \(\alpha_K\), the energy \(e_K\) is zero up to and including \(\alpha_K\) and strictly positive above it. For \(0\le\lambda<\alpha_K\), we can choose \(b\) strictly between \(\lambda\) and \(\alpha_K\). Then \(e_K(b)=0\), so Proposition~\ref{prop:robust} gives \(\Pp(V_N(\lambda)=0)\to1\).  For $\lambda>\alpha_K$, we have $\E V_N(\lambda)/N\to e_K(\lambda)>0$. The concentration bound therefore implies that, with high probability, every assignment violates at least $\tfrac12 e_K(\lambda)N$ clauses.
\end{proof}

Carenini's formulation \cite{Carenini} uses exactly
$m=\lfloor\lambda N\rfloor$ clauses with distinct variables within
each clause. To pass to this convention, we first fix the clause count.
Below the threshold, a Poisson formula at a slightly higher density
still below $\alpha_K$ is satisfiable and has at least $m$ clauses
with high probability, so its first $m$ clauses are satisfiable.
Above the threshold, a Poisson formula at a slightly lower density
still above $\alpha_K$ has at most $m$ clauses with high probability,
and adding the remaining clauses preserves the linear lower bound on
violations. With $m$ fixed, the probability that every clause has
distinct indices is
\[
 \left(\prod_{i=0}^{K-1}\left(1-\frac{i}{N}\right)\right)^m
 \longrightarrow e^{-\lambda\binom{K}{2}}>0.
\]
Conditioning on this event gives exactly Carenini's clause law and
preserves both high-probability conclusions. See also the transfer
for $K=3$ in \cite[proof of Theorem~4.4]{OAI3}.


\begin{thebibliography}{99}

\bibitem{AM}
E. Abbe and A. Montanari,
\emph{On the concentration of the number of solutions of random
satisfiability formulas}, Random Structures Algorithms \textbf{45}
(2014), no.~3, 362--382.
\href{https://arxiv.org/abs/1006.3786}{arXiv:1006.3786}.

\bibitem{AC}
A. Auffinger and W.-K. Chen, W.-K., \emph{Parisi formula for the ground state energy in the mixed p-spin model}, Ann. Probab. \textbf{45} (2017), no.~6B, 4617--4631.
\href{https://arxiv.org/abs/1606.05335}{arXiv:1606.05335}.

\bibitem{BGT}
M. Bayati, D. Gamarnik and P. Tetali,
\emph{Combinatorial approach to the interpolation method and scaling
limits in sparse random graphs}, Ann. Probab. \textbf{41} (2013),
no.~6, 4080--4115.
\href{https://doi.org/10.1214/12-AOP816}{doi:10.1214/12-AOP816}.

\bibitem{BBCKW}
B. Bollob\'as, C. Borgs, J. T. Chayes, J. H. Kim and D. B. Wilson,
\emph{The scaling window of the 2-SAT transition},
Random Structures Algorithms \textbf{18} (2001), 201--256.
\href{https://arxiv.org/abs/math/9909031}{arXiv:math/9909031}.

\bibitem{CareniniWindow}
G. Carenini, \emph{The scaling window of random $K$-SAT},
Electronic Colloquium on Computational Complexity, Report TR26-145 (2026).
\href{https://eccc.weizmann.ac.il/report/2026/145/}{ECCC TR26-145}.

\bibitem{Carenini}
G. Carenini, \emph{A polynomial scaling window for random $K$-SAT and
a proof of the satisfiability conjecture}, Electronic Colloquium on
Computational Complexity, Report TR26-229 (2026).
\href{https://eccc.weizmann.ac.il/report/2026/229/}{ECCC TR26-229}.

\bibitem{CR}
V. Chv\'atal and B. Reed,
\emph{Mick gets some (the odds are on his side)},
Proc. 33rd Annual IEEE Symposium on Foundations of Computer Science
(1992), 620--627.
\href{https://doi.org/10.1109/SFCS.1992.267789}{doi:10.1109/SFCS.1992.267789}.

\bibitem{Coja2}
A. Coja-Oghlan and K. Panagiotou,
\emph{The asymptotic $K$-SAT threshold}, Adv. Math. \textbf{288}
(2016), 985--1068.
\href{https://arxiv.org/abs/1310.2728}{arXiv:1310.2728}.

\bibitem{CGHS}
D. Coppersmith, D. Gamarnik, M. Hajiaghayi and G. B. Sorkin,
\emph{Random MAX SAT, random MAX CUT, and their phase transitions},
Random Structures Algorithms \textbf{24} (2004), 502--545.
\href{https://doi.org/10.1002/rsa.20015}{doi:10.1002/rsa.20015}.

\bibitem{DSS}
J. Ding, A. Sly and N. Sun,
\emph{Proof of the satisfiability conjecture for large $k$},
Ann. of Math. (2) \textbf{196} (2022), no.~1, 1--388.
\href{https://doi.org/10.4007/annals.2022.196.1.1}{doi:10.4007/annals.2022.196.1.1}.

\bibitem{FL}
S. Franz and M. Leone,
\emph{Replica bounds for optimization problems and diluted spin systems},
J. Statist. Phys. \textbf{111} (2003), 535--564.
\href{https://arxiv.org/abs/cond-mat/0208280}{arXiv:cond-mat/0208280}.

\bibitem{Friedgut}
E. Friedgut, \emph{Sharp thresholds of graph properties, and the
$K$-SAT problem}, with an appendix by J. Bourgain,
J. Amer. Math. Soc. \textbf{12} (1999), no.~4, 1017--1054.
\href{https://doi.org/10.1090/S0894-0347-99-00305-7}{doi:10.1090/S0894-0347-99-00305-7}.

\bibitem{Gamarnik}
D. Gamarnik,
\emph{The satisfiability conjecture for random $K$-SAT formulas.
Human-friendly exposition of the recent OpenAI resolution},
preprint (2026).
\href{https://arxiv.org/abs/2611.10139}{arXiv:2611.10139}.

\bibitem{Goerdt}
A. Goerdt, \emph{A threshold for unsatisfiability},
J. Comput. System Sci. \textbf{53} (1996), no.~3, 469--486.
\href{https://doi.org/10.1006/jcss.1996.0081}{doi:10.1006/jcss.1996.0081}.

\bibitem{Guerra}
F. Guerra, \emph{Broken replica symmetry bounds in the mean field spin glass model},
Comm. Math. Phys. \textbf{233} (2003), no.~1, 1--12.
\href{https://doi.org/10.1007/s00220-002-0773-5}{doi:10.1007/s00220-002-0773-5}.

\bibitem{KMRZ}
F. Krzakala, A. Montanari, F. Ricci-Tersenghi, G. Semerjian
and L. Zdeborov\'a,
\emph{Gibbs states and the set of solutions of random constraint
satisfaction problems}, Proc. Natl. Acad. Sci. USA \textbf{104}
(2007), no.~25, 10318--10323.
\href{https://doi.org/10.1073/pnas.0703685104}{doi:10.1073/pnas.0703685104}.

\bibitem{MMZ}
S. Mertens, M. M\'ezard and R. Zecchina,
\emph{Threshold values of random $K$-SAT from the cavity method},
Random Structures Algorithms \textbf{28} (2006), no.~3, 340--373.
\href{https://doi.org/10.1002/rsa.20090}{doi:10.1002/rsa.20090}.

\bibitem{MP1}
M. M\'ezard and G. Parisi,
\emph{The Bethe lattice spin glass revisited},
Eur. Phys. J. B \textbf{20} (2001), 217--233.

\bibitem{MPzero}
M. M\'ezard and G. Parisi,
\emph{The cavity method at zero temperature},
J. Statist. Phys. \textbf{111} (2003), 1--34.
\href{https://arxiv.org/abs/cond-mat/0207121}{arXiv:cond-mat/0207121}.

\bibitem{MPZ}
M. M\'ezard, G. Parisi and R. Zecchina,
\emph{Analytic and algorithmic solution of random satisfiability problems},
Science \textbf{297} (2002), 812--815.

\bibitem{MZcavity}
M. M\'ezard and R. Zecchina,
\emph{Random $K$-SATisfiability problem: From an analytic solution
to an efficient algorithm}, Phys. Rev. E \textbf{66} (2002), 056126.

\bibitem{MPRT}
A. Montanari, G. Parisi and F. Ricci-Tersenghi,
\emph{Instability of one-step replica-symmetry-broken phase in
satisfiability problems}, J. Phys. A \textbf{37} (2004), 2073--2091.
\href{https://arxiv.org/abs/cond-mat/0308147}{arXiv:cond-mat/0308147}.

\bibitem{MRS}
A. Montanari, F. Ricci-Tersenghi and G. Semerjian,
\emph{Clusters of solutions and replica symmetry breaking in random
$K$-SATisfiability}, J. Stat. Mech. (2008), P04004.
\href{https://arxiv.org/abs/0802.3627}{arXiv:0802.3627}.

\bibitem{OAImp}
OpenAI, \emph{The M\'ezard--Parisi formula for diluted spin glasses},
preprint, September 23, 2026.
\href{https://github.com/openai/math/tree/main/preprints/The-Mezard-Parisi-formula-for-diluted-spin-glasses-September-23-2026}{OpenAI Math Release}.

\bibitem{OAIall}
OpenAI, \emph{A limiting satisfiability threshold for every fixed clause
size}, preprint, September 25, 2026.
\href{https://github.com/openai/math/tree/main/preprints/A-Limiting-Satisfiability-Threshold-for-Every-Fixed-Clause-Size-September-25-2026}{OpenAI Math Release}.

\bibitem{OAI3}
OpenAI, \emph{Computing the random 3-SAT threshold}, preprint,
September 27, 2026.
\href{https://github.com/openai/math/tree/main/preprints/Computing-the-Random-3-SAT-Threshold-September-27-2026}{OpenAI Math Release}.

\bibitem{OAIvariance}
OpenAI, \emph{Variance of the random $K$-SAT hitting time}, preprint,
September 27, 2026.
\href{https://github.com/openai/math/tree/main/preprints/Variance-of-the-Random-k-SAT-Hitting-Time-September-27-2026}{OpenAI Math Release}.

\bibitem{SKbook}
D. Panchenko, \emph{The Sherrington--Kirkpatrick Model},
Springer Monographs in Mathematics, Springer, New York (2013).

\bibitem{Panchenko}
D. Panchenko,
\emph{The M\'ezard--Parisi formula}, preprint (2026).
\href{https://arxiv.org/abs/2610.39545}{arXiv:2610.39545}.

\bibitem{PT}
D. Panchenko and M. Talagrand,
\emph{Bounds for diluted mean-fields spin glass models},
Probab. Theory Relat. Fields \textbf{130} (2004), no.~3, 319--336.
\href{https://arxiv.org/abs/math/0405357}{arXiv:math/0405357}.

\bibitem{Talagrand}
M. Talagrand, \emph{Mean Field Models for Spin Glasses, Volume I:
Basic Examples}, Ergebnisse der Mathematik und ihrer Grenzgebiete,
vol.~54, Springer, Berlin (2011).

\end{thebibliography}
\end{document}